\documentclass[12pt]{article}
\usepackage{setspace}
\usepackage{amsfonts}
\usepackage{amsmath}
\usepackage{amssymb}
\usepackage{amsthm}
\usepackage{mathrsfs}
\usepackage[all]{xy}

\newcommand{\proofend}{\hfill \hbox{\vrule width 5pt height 5pt depth
0pt}}
\newcommand{\R}{\mathbb{R}}
\newcommand{\C}{\mathbb{C}}
\newcommand{\Z}{\mathbb{Z}}
\newcommand{\N}{\mathbb{N}}
\newcommand{\Q}{\mathbb{Q}}

\newcommand{\proj}{\mathbb{P}}
\newcommand{\A}{\mathbb{A}}

\newcommand{\spec}{\mathrm{Spec} \,}

\begin{document}

\begin{center}

{\large \sc A lower bound on the orbit growth \\ of a regular self-map of affine space}

\vspace{1cm}

Vesselin Dimitrov \\

\medskip

\texttt{vesselin.dimitrov@yale.edu}

\bigskip

\vspace{1cm}

\end{center}

\begin{abstract}
We show that if $f : \A_{\bar{\Q}}^r \to \A_{\bar{\Q}}^r$ is a regular self-map and $P \in \A^r(\bar{\Q})$ has $\limsup_{n \in \N} \frac{\log{h_{\mathrm{aff}}(f^nP)}}{\log{n}}
< 1/r$, where $h_{\textrm{aff}}$ is the affine logarithmic Weil height, then $\N$ partitions into a finite set and finitely many full arithmetic progressions, on each of which the coordinates of $f^nP$ are polynomials in $n$.

In particular, if $(f^nP)_{n \in \N}$ is a Zariski-dense orbit, then either $r = 1$ and $f$ is of the shape $t \mapsto \zeta t + c$, $\zeta \in \mu_{\infty}$, or else $\limsup_{n \in \N} \frac{\log{h_{\mathrm{aff}}(f^nP)}}{\log{n}} \geq 1/r$. This inequality is the exponential improvement of the trivial lower bound obtained from counting the points of bounded height in $\A^r(K)$.
\end{abstract}

\vspace{1cm}

\begin{center}
{\bf 1. Introduction}
\end{center}

{\bf 1.1. } In the appendix to our preprint~\cite{prelim} we formulated a precise conjectural criterion for the algebraicity of a formal function on a projective curve over a global field. For the case of the projective line, this criterion generalizes simultaneously the classical P\'olya-Bertrandias criterion (cf. Amice~\cite{amice}, Ch.~5), on the one hand, and a conjecture of I. Ruzsa~\cite{rusza,perellizannier,zannierruzsa,christol}, on the other hand. In~\cite{prelim} (3.2 in~\emph{loc. cit.}) we proved a weak variant of this conjecture with the purpose of applying it to a case of a generalization of the ``Hadamard quotient theorem'' to higher genus and positive equicharacteristic (Theorem~1.7 in {\it loc. cit.}). In the present note we prove another weak variant, and apply it to obtain a lower bound on the growth of an orbit under a regular iteration of affine space.

\medskip

{\bf 1.2. } We need to settle some notation before we can state our result. In what follows we denote by $h(\cdot) : \bar{\Q}^{\times} \to \R^{\geq 0}$ the absolute logarithmic Weil height, whose definition we recall next. For $p$ a finite rational prime let $|\cdot|_p$ by the $p$-adic absolute value on $\Q$ normalized by $|p|_{p} = 1/p$, while for $p = \infty$ we consider the ordinary (archimedean) absolute value. For $v$ a place a number field $K$ lying over the place $p$ of $\Q$, and for $x \in K$, let $|x|_{v} := |N_{K_v/\Q_p}(x)|_p^{1/[K:\Q]}$.
 For $\mathbf{x} = [x_0:x_1:\cdots:x_r] \in \proj^r(K)$ we then have
\begin{eqnarray*}
h(\mathbf{x}) := \sum_v \max_j \log{|x_j|_v}.
\end{eqnarray*}
This definition is independent of the choice of the number field $K$.
Viewing $\A^r(\bar{\Q})$ as the affine piece $[1:x_1:\ldots:x_r]$ in $\proj^r(\bar{\Q})$, we consider the affine height $h_{\mathrm{aff}}(\alpha_1,\ldots,\alpha_r) := h([1:\alpha_1:\cdots:\alpha_r])$. Finally, we write as usual $h(\alpha) = h_{\mathrm{aff}}(\alpha)$ for $\alpha \in \bar{\Q} = \A^1(\bar{\Q})$.

For a polynomial $F \in K[\mathbf{x}]$ in several variables over $K$, we write $h(F)$ for the height of its set of coefficients, viewed as a point in a projective space. For $F = \sum_{n \geq 0} a_nt^n \in K[[t]]$, let $F_{/N} := \sum_{n=0}^N a_nt^n \in K[t]$ be the polynomial truncation modulo $t^{N+1}$.

\medskip

{\bf 1.3. } To motivate our result we make the following trivial observation. If a point $P \in \Z^r$ has infinite orbit under a set-theoretic mapping $f: \Z^r \to \Z^r$, then $H_n := \max_{0 \leq i \leq n} \exp(h(f^iP))$ satisfies $(2H_n+1)^r > n$, whence, in the limit,
$$
\limsup_{n \in \N} \frac{h_{\mathrm{aff}}(f^nP)}{\log{n}} \geq 1/r.
$$
The result of this note is that when the mapping $f$ has an algebraic structure, and outside of a degenerate situation, this trivial inequality can be improved exponentially.

\medskip

{\bf Theorem 1.4. }  {\it Consider regular maps $f : \A_{\bar{\Q}}^r \to \A_{\bar{\Q}}^r$ and $\lambda : \A_{\bar{\Q}}^r \to \A_{\bar{\Q}}^1$. Let $P \in \A^r(\bar{\Q})$ be a point such that
$$
\limsup_{n \in \N} \frac{\log{h(\lambda(f^nP))}}{\log{n}} < 1/r.
$$
Then there is a $d \in \N$ and polynomials $p_0,\ldots,p_{d-1} \in \bar{\Q}[x]$ such that $\lambda(f^{nd+i}P) = p_i(n)$ for all $i = 0, \ldots,d-1$ and $n \gg 0$.}

The following is an immediate corollary, obtained by taking $\lambda$ to be the coordinate projections.

\medskip

{\bf Corollary 1.5. } {\it In the setup of Theorem~1.3, either the Zariski closure of the orbit $(f^nP)_{n \in \N}$ is a union of rational curves in $\A^r$, or else
$$
\limsup_{n \in \N} \frac{\log{h_{\mathrm{aff}}(f^nP)}}{\log{n}} \geq 1/r.
$$}

\medskip

Another corollary arises in taking $f$ to be a recurrence of the form
$$
(x_1,\ldots,x_{r-1},x_r) \mapsto (x_2,\ldots,x_r,p(x_1,\ldots,x_r))
$$
 and $\lambda$ the projection onto the last coordinate.

\medskip

{\bf Corollary 1.6. } {\it Consider a sequence $A : \N \to \Z$ satisfying a finite polynomial recurrence $A(n+r) = p(A(n+r-1),\ldots,A(n))$, where $p \in \Z[T_1,\ldots,T_r]$. Then either
$$
\limsup_{n \in \N} \frac{\log\log{|A(n)|}}{\log{n}} \geq 1/r,
$$
or else there is a $d \in \N$ such that each of the sequences $(A(nd+i))_{n \geq 0}$, $0 \leq i < d$, is a polynomial for $n \gg 0$.}

\medskip

We will derive Theorem~1.4 from a criterion for the rationality of a formal function on the projective line, which we formulate and prove in the next section. Results of this type date back to the short paper~\cite{zannierposk} by U. Zannier.

\bigskip

\begin{center}
{\bf 2. A rationality criterion}
\end{center}

In what follows, $K$ is a number field and $S$ a finite set of places of $K$ including all archimedean places. The following is a particular case of the conjecture formulated in the appendix to~\cite{prelim}.

\medskip

{\bf Conjecture 2.1. } {\it Consider $F \in O_{K,S}[[t]]$ a formal power series with $S$-integral coefficients which has a positive radius of convergence $R_v > 0$ at each place $v$ of $K$. For each prime (closed point) $s \in  \spec{O_K}$ let $h_s \in \N_0 \cup \{\infty\}$ be $+\infty$ if either $s \in S$ or the reduction $F \mod{s} \in k(s)[[t]]$ is not in $k(s)(t)$; and the degree of the rational function $F\mod{s}$ otherwise. If
\begin{equation} \label{shp}
\sum_{v} \log{R_v} + \liminf_{n \in \N} \Big\{ \frac{1}{n} \frac{1}{[K:\Q]} \sum_{s : \,h_s < n} \log{|k(s)|} \Big\} > 0,
\end{equation}
then the power series $F \in K(t)$ is rational.}

 \medskip

This conjecture is sharp, in the sense that there are uncountably many power series $F$ for which the quantity on the left-hand side of~(\ref{shp}) is zero. It extends an old conjecture of I. Ruzsa~\cite{rusza}, and appears intractable to our current techniques. In this section we prove the following crude variant.

\medskip

{\bf Proposition~2.2. } {\it Let $\eta > 0$. In the setup of Conjecture~2.1, assume instead that the inequality
\begin{equation} \label{strongeq}
  \frac{1}{[K:\Q]}  \sum_{s \, : \, h_s < n/(2+\eta)} \log{|k(s)|}  > \frac{3}{2}\log{n} + \Big( 1+\frac{1}{\eta} \Big) h(F_{/n}) + \frac{1}{2}\log{|D_{K/\Q}|}
\end{equation}
holds for all $n \gg 0$.
Then $F \in K(t)$ is rational. }

\medskip

{\it Proof. } The proof is a variant of that presented to the algebraicity criterion~3.2 in~\cite{prelim}. Without loss of generality we may assume $\eta$ to be rational. Letting $L$ be a large integer parameter such that $\eta L \in \Z$, Siegel's lemma (see 3.1 in~\cite{prelim} and the references therein) for $M := L$ equations in $N := (1+\eta) L$ unknowns yields polynomials $P \in K[t]$ and $Q \in K[t] \setminus \{0\}$ of degrees less than $(1+\eta)L$ such that
\begin{equation}
Q(t)F(t) - P(t) \equiv 0 \mod{t^{ (2+\eta)L }}
\end{equation}
and
\begin{eqnarray} \label{estimate}
h(Q) \leq \frac{1}{2}\log{N} + \frac{1}{2}\log{|D_{K/\Q}|} + \frac{1}{\eta} h\big( F_{/(2+\eta)L} \big).
\end{eqnarray}
 It follows from~(\ref{strongeq}) and (\ref{estimate}) that there is an $L < \infty$ such that all $n \geq (2+\eta)L$ satisfy
 \begin{equation} \label{contr}
\frac{1}{[K:\Q]} \sum_{s : \, h_s < n/(2+\eta)} \log{|k(s)|} > \log{N} + h(Q) + h(F_{/n}).
 \end{equation}
 {\it We claim that then $F = P/Q$, hence $F$ is rational.} Assuming otherwise, let $n$ be the minimum integer such that $Q(t)F(t) - P(t) \not\equiv 0 \mod{t^n}$; by construction, $n > (2+\eta)L$. Consider a prime $s$ of $K$ with $h_s < n/(2+\eta)$, and write $F_s := A_s/B_s$
  the reduction at $s$, with $A_s, B_s \in k(s)[t]$, $\deg{A_s}, \deg{B_s} < n/(2+\eta)$. Then, denoting by a tilde the reduction at $s$, we have $A_s(t)\widetilde{Q}(t) - B_s(t)\widetilde{P}(t) \equiv 0 \mod{t^{n-1}}$. The degree of this polynomial is less than $(1+\eta)L + n/(2+\eta)$, which by our assumption $n > (2+\eta)L$ does not exceed $n-1$. It follows that the polynomial is identically zero, hence the coefficients of $Q(t)F(t) - P(t)$ all reduce to zero at $s$.

  On the other hand, we have assumed that $t^n$ appears with a non-zero coefficient $c \in K \setminus \{0\}$
 in $Q(t)F(t) - P(t)$. (This is just the coefficient of $t^n$ in $Q(t)F(t)$.) Thus the product formula yields
 \begin{equation} \label{productform}
 \sum_{\textrm{all } v} -\log{|c|_v} = 0.
\end{equation}
 At the places corresponding to the primes $s$ with $h_s < n/(2+\eta)$, the previous paragraph shows that the contribution to~(\ref{productform}) is at least $\frac{1}{[K:\Q]}\log{|k(s)|}$. The absolute value of the sum of the remaining contributions does not exceed $h(c)$,  and~(\ref{productform}) yields the lower bound
 \begin{equation}  \label{lower}
 h(c) \geq \frac{1}{[K:\Q]} \sum_{s : \, h_s < n/(2+\eta)} \log{|k(s)|}.
 \end{equation}
To estimate $h(c)$ from above, we use the easy bound
$$
h(\alpha_1 + \cdots + \alpha_r) \leq \log{r} + \sum_v \max_j \log^+{|\alpha_j|_v},
$$
applied to the sum defining $c$ as the coefficient of $t^n$ in $Q(t)F(t)$. We obtain:
\begin{equation} \label{upper}
h(c) \leq \log{N} + h(Q) + h(F_{/n}).
\end{equation}
Taken together with~(\ref{lower}) this contradicts~(\ref{contr}), thus forcing $F = P/Q$ as claimed. \proofend

 \bigskip

 \begin{center}
 {\bf 3. Proof of Theorem~1.4 }
 \end{center}

 There is a number field $K$ and a finite set $S$ of its places, including all archimedean places, such that the triple $(f,\lambda,P)$ has a model over $O_{K,S}$.
 We will apply Proposition~2.2 to the formal power series $\Phi := \sum_{n \geq 0} \lambda(f^nP) t^n \in O_{K,S}[[t]]$. If the power series $\Phi$ is rational, the conclusion of Theorem~1.4 follows from the explicit descriptions of coefficients of rational power series as confluent power sums. If $\Phi$ is not rational, Proposition~2.2 with $\eta := 1$ implies the lower bound inequality
 \begin{equation} \label{reverse}
 h(\Phi_{/n}) \geq  \frac{1}{2[K:\Q]}  \sum_{s \, : \, h_s(\Phi) < n/3} \log{|k(s)|} - O(\log{n})
 \end{equation}
 for infinitely many $n \in \N$.
 On the other hand, for $s \notin S$ a prime of $K$, the iteration $f : \A_{O_{K,S}}^r \to \A_{O_{K,S}}^r$ reduces mod $s$ to an iteration $f\mod{s} : \A_{k(s)}^r \to \A_{k(s)}^r$ over a set with $|k(s)|^r$ elements, hence $h_s(\Phi) \leq 2|k(s)|^r$. Consequently, by the prime number theorem, (\ref{reverse}) yields
 \begin{eqnarray} \label{ataka}
h(\Phi_{/n}) \geq  \frac{1}{3}  (n/6)^{1/r} \quad \textrm{for arbitrarily large } n \in \N.
 \end{eqnarray}
 We have $h(\Phi_{/n}) \leq |S|\max_{j \leq n} h(\lambda(f^jP))$, and the conclusion of the theorem follows from~(\ref{ataka}). \proofend

\bigskip

\begin{center}
{\bf 4. Two conjectures}
\end{center}

 We end this note by recording two conjectures related to the setup of Theorem~1.4. Problems of this type have been posed by J.H. Silverman in~\cite{silvrational}.

 \medskip

 {\bf Conjecture~4.1. } {\it Let $X$ be a complex projective variety, $f : X \dashrightarrow X$ a rational self-map, $\lambda : X \dashrightarrow \proj^1$ a non-constant rational function, and $P \in X(\C)$ a point with well-defined forward orbit. Then the set $\{ n \mid \lambda(f^nP) = 0 \} \subset \N_0$ is the union of a finite set with finitely many full arithmetic progressions. }

 \medskip

 {\bf Conjecture~4.2. } {\it Let $A/\bar{\Q}$ be an abelian variety  and $\lambda : A \dashrightarrow \proj^1$ a non-constant rational function. Consider a point $P \in A(\bar{\Q})$. If $h\big( \lambda([n]P) \big) = o(n^2)$ then there is a surjective homomorphism $A \twoheadrightarrow B$ to an abelian variety mapping $P$ to a torsion point.}

\end{document}